\documentclass[fleqn,11pt,oneside]{article}

\usepackage{amsmath, amsthm, amssymb}
\usepackage{pdfpages}
\usepackage{xcolor}
\usepackage{graphicx}

\newcommand{\sa}{\sigma}
\newcommand{\la}{\lambda}
\newcommand{\eps}{\varepsilon}
\newcommand{\R}{\mathbb{R}}
\newcommand{\C}{\mathbb{C}}
\newcommand{\mtx}{}

\newcommand{\wt}{\widetilde}

\newcounter{algo}

\newtheorem{algorithm}[algo]{Algorithm}
\newtheorem{theorem}{Theorem}[section]
\newtheorem{lemma}[theorem]{Lemma}
\newtheorem{remark}[theorem]{Remark}
\newtheorem{observe}[theorem]{Observation}

\begin{document}

\thispagestyle{empty}
\begin{center}
  \begin{minipage}[t]{6.0in}

We introduce an efficient scheme for the construction of quadrature rules
for bandlimited functions.
While the scheme is predominantly based on well-known facts about prolate
spheroidal wave functions of order zero, it has the asymptotic CPU time
estimate $O(n \log n)$ to construct an $n$-point quadrature rule. 
Moreover, the size of the ``$n \log n$'' term in the CPU time estimate is
small, so for all practical purposes the CPU time cost is proportional to
$n$. 
The performance of the algorithm is illustrated by several numerical
examples.

 \vspace{ -100.0in}

 \end{minipage}

\end{center}

\vspace{ 4.60in}

\begin{center}

  \begin{minipage}[t]{4.4in}

\begin{center}

{\bf
A fast procedure for the construction of quadrature formulas for bandlimited
functions
} \\

 \vspace{ 0.50in}

A. Gopal$\mbox{}^{\dagger, \star}$,  V. Rokhlin$\mbox{}^{\dagger, \ddagger}$ \\
                May 2, 2023
\end{center}

 \vspace{ -100.0in}

 \end{minipage}

          \end{center}

\vfill

\noindent $\mbox{}^{\star}$ This author was supported in part by AFOSR
(grant no.~FA9550-19-1-0320).

\noindent $\mbox{}^{\ddagger}$ This author was supported in part by NSF (grant no.~DMS-1952751) and ONR (grant no.~N00014-14-1-0797).

\vspace{2mm}

\noindent
$\mbox{}^{\dagger}$ Dept.~of Mathematics, Yale University, New Haven CT 06511

\vspace{2mm}

\noindent
Approved for public release: distribution is unlimited.


\pagebreak
\clearpage
\pagenumbering{arabic}

\section{Introduction}
A function $f: \R \to \C$ is said to be \emph{bandlimited} with
\emph{bandlimit} $c > 0$ if there exists $\sigma \in L^2([-1,1])$ such that
\begin{equation}
f(x) = \int_{-1}^1 e^{ictx} \sa(t)\,dt, \quad x \in \R.
\label{eq:bandlimfun}
\end{equation}
In other words, a function is bandlimited if its Fourier transform is
compactly supported. 
Such functions arise in a variety of applications such as fluid dynamics,
signal processing, and scattering theory. 
Typically, only the restriction of $f$ to some compact domain is of
interest, such as an interval.

In this regime, the standard techniques of Chebyshev and Legendre
interpolation and quadrature can often be effectively employed. 
However, polynomials are a suboptimal basis for the space of
bandlimited functions, and thus, schemes based on polynomial approximation
are inherently suboptimal in the number of quadrature points
required to obtain a prescribed accuracy (cf., \cite{boyd2004prolate}). 
Indeed, the optimal basis for bandlimited functions on an interval is
actually the prolate spheroidal wavefunctions of order zero (PSWFs).
The PSWFs were introduced by Slepian, Landau, and Pollack in the context of
information theory in
\cite{slepian1961prolate,landau1961prolate,landau1962prolate,slepian1964prolate,slepian1978prolate},
and more recently there have been several advances in numerical algorithms
for PSWFs (e.g.,
\cite{xiao2001prolate,shkolnisky2006approximation,osipov2013numerical}). 

In this paper, we continue this latter line of work and present an
efficient and robust numerical algorithm for computing quadrature rules for
bandlimited functions on an interval.
For any fixed bandlimit $c > 0$, the algorithm computes an $n$-point
quadrature rule in $O(n \log n)$ floating point operations,
unlike the generalized Gaussian quadrature schemes of
\cite{xiao2001prolate} and \cite{shkolnisky2006approximation} which require
$O(n^3)$ operations and the scheme of \cite{osipov2013numerical} which
requires $O(n^2)$ operations. 
The nodes of the computed rules are taken to be the roots of appropriately
chosen PSWFs, and the weights are chosen to exactly integrate certain
collections of PSWFs.
The rules are stable and require substantially fewer nodes to integrate
bandlimited functions to a prescribed accuracy when compared to classical
polynomial quadrature.


This paper is organized as follows.
Section \ref{s:prelims} summarizes some standard facts needed in this paper.
Section \ref{s:apparatus} contains various mathematical results and
numerical techniques used in our algorithm.
We describe our algorithm in Section \ref{s:algorithm},
and the results of several numerical experiments are provided in Section
\ref{s:numerics}.
Section \ref{s:conclusions} contains our conclusions. 

\section{Preliminaries}
\label{s:prelims}

\subsection{Legendre polynomials and functions}
\label{s:legendre}

The Legendre polynomials $\{P_k\}_{k=0}^\infty$ are defined by the
three-term recurrence relation
\begin{equation}
P_{k+1}(x) = \frac{2k+1}{k+1} x P_k(x) - \frac{k}{k+1} P_{k-1}(x), \quad k = 1,2,\hdots,
\label{eq:legerec}
\end{equation}
with the initial conditions
\begin{equation}
P_0(x) = 1, \quad P_1(x) = x.
\end{equation}
Differentiating \eqref{eq:legerec} yields the recurrence relation
\begin{equation}
P_k'(x) = \frac{2k-1}{k} P_{k-1}(x) + \frac{2k-1}{k} x P_{k-1}'(x) -
\frac{k-1}{k} P_{k-2}'(x), \quad k = 2,3,\hdots.
\label{eq:diffrec}
\end{equation}

The Legendre functions of the second kind $\{Q_k\}_{k=0}^\infty$ are also
defined by the recurrence relation
\begin{equation}
Q_{k+1}(x) = \frac{2k+1}{k+1} x Q_k(x) - \frac{k}{k+1} Q_{k-1}(x), \quad k = 1,2,\hdots,
\end{equation}
but with the initial conditions
\begin{equation}
Q_0(x) = \frac{1}{2} \log \frac{1+x}{1-x}, \quad Q_1(x) = \frac{x}{2} \log
\frac{1+x}{1-x}-1.
\end{equation}
For any $z \in \C$ and integer $k \geq 0$, the identity
\begin{equation}
Q_k(z) = \frac{1}{2} \int_{-1}^1 \frac{P_k(t)}{z-t}\,dt
\label{eq:legeid}
\end{equation}
holds (\cite[8.8.3]{abramowitz1988handbook}).

\begin{remark}
For any integer $k \geq 0$, $P_k$ and $Q_k$ are independent solutions to the
ODE 
\begin{equation}
(1-x^2)y''(x)-2xy'(x)+ k(k+1)y(x)=0, \quad x \in [-1,1].
\end{equation}
\end{remark}

\subsection{Pr\"ufer transform}
Suppose that the function $\psi$ satisfies a second-order ordinary
differential equation of the form 
\begin{equation}
p(x) \psi''(x) + q(x) \psi'(x) + r(x) \psi(x) = 0, \quad x \in (-1,1),
\label{eq:gen_ode}
\end{equation}
where $p$, $q$, and $r$ are polynomials of degree 2 and $p$ and $r$ are
positive.
Given a differentiable, positive function $\gamma : \R \to \R$, we
can define the function $\theta : \R \to \R$ by the formula
\begin{equation}
\theta(x) = \arctan \left( \frac{1}{\gamma(x)} \frac{p(x)
\psi'(x)}{\psi(x)}\right) + m \pi,
\label{eq:gen_theta}
\end{equation}
where $m$ is an arbitrary integer.
Substituting \eqref{eq:gen_theta} into \eqref{eq:gen_ode} produces the
first-order ODE
\begin{equation}
\theta' = -\frac{\gamma}{p} \sin^2(\theta) -
\frac{r}{\gamma} \cos^2(\theta) - \left(
\frac{\gamma'}{\gamma} + \frac{q-p'}{p} \right) \frac{\sin(2 \theta)}{2},
\quad x \in (-1,1).
\label{eq:gen_prufode}
\end{equation}
From \eqref{eq:gen_theta}, it immediately follows that the roots of $\psi$
are the values of $x \in (-1,1)$ for which $\theta(x) = (k+1/2)\pi$ for some
integer $k$; 
similarly, the roots of $\psi'$ are the values of $x$ for which $\theta(x) =
k \pi$ for some integer $k$.

It is convenient to choose 
\begin{equation}
\gamma(x) = \sqrt{r(x)p(x)}.
\label{eq:gamma}
\end{equation}
Then, \eqref{eq:gen_prufode} simplifies to 
\begin{equation}
\theta' = -\sqrt{\frac{r}{p}} - \frac{r'p-p' r + 2rq}{2rp} \frac{\sin(2
\theta)}{2}.
\label{eq:gen_prufodereal}
\end{equation}
Noting in \eqref{eq:gen_prufodereal} that $d\theta/dx < 0$ for all $x \in
(-1,1)$, we take the reciprocal, which yields
\begin{equation}
\frac{dx}{d\theta} =  
- \left(\sqrt{\frac{r}{p}} + \frac{r'p-p' r + 2rq}{2rp} \frac{\sin(2
\theta)}{2}\right)^{-1}.
\label{eq:gen_prufode2}
\end{equation}

\begin{observe}
Given $x_0$, a root of $\psi$, the ODE \eqref{eq:gen_prufode2} can be used
to determine all other roots by solving \eqref{eq:gen_prufode2} subject to
the initial condition
\begin{equation}
x \left( \frac{\pi}{2} \right) = x_0
\end{equation}
and evaluating the solution at half-integer multiples of $\pi$. 
Alternatively, if $x_0'$, a root of $\psi'$, is available, the initial
condition 
\begin{equation}
x \left( 0 \right) = x_0'
\end{equation}
can be used instead.
\end{observe}
 

\subsection{Taylor series method for ODEs}
\label{s:taylor}

For any integer $k \geq 0$, differentiating \eqref{eq:gen_ode} $k$ times
produces the recurrence relation
\begin{align}
p\psi^{(k+2)} = &-(kp'+q) \psi^{(k+1)} - \left( \frac{k(k-1)}{2} p'' + kq' + r \right) \psi^{(k)} \nonumber \\
&-\left( \frac{k(k-1)}{2} q'' + k r' \right) \psi^{(k-1)} - \frac{k(k-1)}{2} r'' \psi^{(k-2)},
\label{eq:gen_derivs}
\end{align}
where we define $\psi^{(-2)}(x)=0$ and $\psi^{(-1)}(x)=0$.
Given a point $x_* \in (-1,1)$ and the values of $\psi(x_*)$ and
$\psi'(x_*)$, higher-order derivatives of $\psi$ at $x_*$ can be evaluated
using \eqref{eq:gen_derivs}.
Fixing an integer $P \geq 2$ and a sufficiently small $h > 0$, $\psi$ and
$\psi'$ can then be evaluated at $x_*+h$ using the Taylor series
\begin{equation}
\psi(x_*+h) = \sum_{k=0}^P \frac{\psi^{(k)}(x_*)}{k!} h^k + O(h^{P+1}),
\label{eq:gen_taylor1}
\end{equation}
and
\begin{equation}
\psi'(x_*+h) = \sum_{k=0}^{P-1} \frac{\psi^{(k+1)}(x_*)}{k!} h^k + O(h^{P}).
\label{eq:gen_taylor2}
\end{equation}

This provides a method for solving the initial value problem
\eqref{eq:gen_ode} with the initial conditions
\begin{equation}
\psi(x_0) = \psi_0, \quad \psi'(x_0) = \psi_0'.
\end{equation}
To determine the values of $\psi$ and $\psi'$ at the next time step $x_1$,
the recurrence relation \eqref{eq:gen_derivs} is used to compute $\{
\psi^{(k)}(x_0)\}_{k=2}^P$, and the series \eqref{eq:gen_taylor1} and
\eqref{eq:gen_taylor2} with $x_* = x_1$ and $h=x_1-x_0$ are used to evaluate
$\psi(x_1)$ and $\psi'(x_1)$.
This process is repeated for each time step. 

\subsection{Prolate spheroidal wave functions}
\textbf{Basic properties}

\noindent For a fixed \emph{bandlimit} $c > 0$, prolate spheroidal wave
functions of order zero with bandlimit $c$ (PSWFs) $\{\psi_j\}_{j=0}^\infty$
are eigenfunctions of the truncated Fourier transform, defined by the
formula
\begin{equation}
F_c[\sigma](x) = \int_{-1}^1 e^{i c t x} \sigma(t)\,dt, \quad x \in [-1,1].
\label{eq:fourop}
\end{equation}
The operator $F_c$ in \eqref{eq:fourop} is viewed as acting on
$L^2([-1,1])$.
The PSWFs can be chosen to be real-valued, and we normalize them so that $\|
\psi_j\|_2 = 1$ (here, $\| f \|_2$ denotes the $L^2([-1,1])$-norm of
$f$). 
Moreover, they constitute an orthogonal basis for $L^2([-1,1])$.

Since the PSWFs are analytic on $[-1,1]$, it follows that if $\psi_n$ is
expanded into a Legendre series of the form
\begin{equation}
\psi_n(x) = \sum_{k=0}^\infty \alpha_k P_k(x), \quad x \in [-1,1],
\label{eq:legeexp}
\end{equation}
where $P_k$ is the Legendre polynomial of degree $k$ (see Section
\ref{s:legendre}), then the coefficients $\{\alpha_k\}_{k=0}^\infty$ decay
superalgebraically.
Similarly, $\psi_n'$ can also be expanded into the series
\begin{equation}
\psi_n'(x) = \sum_{k=1}^\infty \alpha_k P_k'(x), \quad x \in [-1,1],
\label{eq:diffexp}
\end{equation}
where $\{P_k'\}_{k=1}^\infty$ can be evaluated using \eqref{eq:diffrec}.

The corresponding eigenvalues, $\{\lambda_j\}_{j=0}^\infty$, are simple, and  
we use the standard ordering convention so that 
\begin{equation}
|\la_0| \geq |\la_1| \geq \hdots \geq 0,
\end{equation}
which results in $\psi_j$ being odd and even for odd and even $j$,
respectively.
For large $c$, there are about $2c/\pi$ eigenvalues with modulus close to
$\sqrt{2\pi/c}$, and the remaining eigenvalues
decay superalgebraically with respect to the index (cf.~\cite[Theorem
2.4]{osipov2013numerical}).

The PSWFs are also the eigenfunctions of a differential operator. 
In particular, the function $\psi_n$ satisfies the equation
\begin{equation} 
(1-x^2) \psi_n''(x) - 2x \psi_n'(x) + (\chi_n - c^2 x^2) \psi_n(x) = 0, \quad x \in (-1,1)
\label{eq:ode}
\end{equation}
for some $\chi_n > 0$.
Note that $\{\chi_j\}_{j=1}^\infty$ is a positive, strictly increasing
sequence.

For a comprehensive overview of the properties of PSWFs, we refer the reader
to \cite{osipov2013numerical}.

\noindent \textbf{Recurrence relations for higher-order derivatives}

\noindent For any integer $k \geq 0$, differentiating \eqref{eq:ode} $k$
times yields the recurrence relation
\begin{align}
(1-x^2) \psi_n^{(k+2)}(x) &= 2x (1 + k) \psi_n^{(k+1)}(x) + (k^2 + k + c^2
x^2 - \chi_n) \psi_n^{(k)}(x) \nonumber \\
&+ 2c^2 k x \psi_n^{(k-1)}(x) + c^2 k (k - 1) \psi_n^{(k-2)}(x),
\label{eq:derivs}
\end{align}
where $\psi_n^{(-1)}(x) = 0$ and $\psi_n^{(-2)}(x) = 0$. 

\begin{observe}
The recurrence relation \eqref{eq:derivs} provides a way of evaluating
$\psi_n$ and $\psi_n'$ in the vicinity of a point $x_*$ where their values
are given, using the approach outlined in Section \ref{s:taylor}.
Specifically, for an integer $P \geq 2$, \eqref{eq:derivs} is used to
evaluate $\{\psi_n^{(k)}(x_*)\}_{k=2}^P$, and then \eqref{eq:gen_taylor1}
and \eqref{eq:gen_taylor2} are used with $\psi=\psi_n$ to evaluate $\psi_n$
and $\psi_n'$ at points near $x_*$.
\end{observe}


\noindent \textbf{Pr\"ufer transform}

\noindent The PSWF equation \eqref{eq:ode} has the Pr\"ufer transform 
\begin{equation}
\frac{d\theta}{dx} = - \sqrt{\frac{\chi_n - c^2 x^2}{1-x^2}} + \frac{c^2 x -
2c^2 x^3 +x \chi_n}{2c^2 x^4 + 2 \chi_n - 2c^2 x^2 - 2x^2 \chi_n} \sin(2
\theta).
\end{equation}
Then,
\begin{equation}
\theta(x) = \arctan \left( \sqrt{\frac{1-x^2}{\chi_n - c^2 x^2}}
\frac{\psi_n'(x)}{\psi_n(x)} \right)+ m \pi
\end{equation}
for an arbitrary integer $m$.
The corresponding analog of \eqref{eq:gen_prufode2} is
\begin{equation}
\frac{dx}{d\theta} = - \left( \sqrt{\frac{\chi_n - c^2 x^2}{1-x^2}} -
\frac{c^2 x - 2c^2 x^3 +x \chi_n}{2c^2 x^4 + 2 \chi_n - 2c^2 x^2 - 2x^2
\chi_n} \sin(2 \theta) \right)^{-1}.
\label{eq:prufode2}
\end{equation}

\subsection{Sturm bisection}
\label{s:sturm}

The following elementary results are the basis of the root-finding technique
known as Sturm bisection and can be found on pages 227--229 of
\cite{trefethen1997numerical}, for example.

\begin{lemma}
Let $T \in \R^{n \times n}$ be a symmetric tridiagonal matrix, and let $T_r
\in \R^{r \times r}$ denote the $r$th leading principal submatrix of $T$.
Moreover, let $p_r$ denote the characteristic polynomial of $T_r$, and define
$p_{-1}(x)=0$ and $p_0(x)=1$.
Then,
\begin{equation}
p_k(x) = (T(r-1,r-1)-x) p_{r-1}(x) - T(r-1,r-2)^2 p_{r-2}(x), \quad r =
2,3,\hdots,n,
\end{equation}
where we have zero-indexed $T$.
\label{lemma:charpol}
\end{lemma}

\begin{theorem}
Let $T \in \R^{n \times n}$ be a symmetric tridiagonal matrix with no zero
subdiagonal entries, and let $T_r \in \R^{r \times r}$ denote the
$r$th leading principal submatrix of $T$.
In addition, let $p_r$ denote the characteristic polynomial of $T_r$ and
define
$p_0(x) = 1$.
Let $\rho(x)$ denote the number of sign changes in the sequence
\begin{equation}
 p_0(x), p_1(x), \hdots, p_n(x)
\end{equation}
with the convention that $p_r(x)$ has the opposite sign of $p_{r-1}(x)$ if
$p_r(x) = 0$. 
Then, $\rho(x)$ equals the number of eigenvalues of $T$ less than $x$.
\label{thm:sturm}
\end{theorem}

\section{Analytical and numerical apparatus}
\label{s:apparatus}

\subsection{Computing PSWFs}

The following theorem states that the coefficients in the Legendre expansion
of $\psi_n$ can be obtained by solving an eigenvalue problem.
We define the infinite, symmetric tridiagonal matrices
\begin{align}
\mtx{A}_{\rm even}(k,k) &= 2k (2k+1) + \frac{4k(2k+1)-1}{(4k+3)(4k-1)} c^2  \label{eq:aeven1} \\
\mtx{A}_{\rm even}(k,k+1) = \mtx{A}_{\rm even}(k+1,k) &= \frac{(2k+2)(2k+1)}{(4k+3)\sqrt{(4k+1)(4k+5)}} c^2  \label{eq:aeven2}, 
\end{align}
and
\begin{align}
\mtx{A}_{\rm odd}(k,k) &= (2k+1)(2k+2) + \frac{(4k+2)(2k+2)-1}{(4k+5)(4k+1)} c^2 \label{eq:aodd1} \\
\mtx{A}_{\rm odd}(k,k+1) = \mtx{A}_{\rm odd}(k+1,k) &= \frac{(2k+3)(2k+2)}{(4k+5)\sqrt{(4k+3)(4k+7)}} c^2\label{eq:aodd2}
\end{align}
for $k=0,1,\hdots$.

The following theorem is the basis for the algorithm in
\cite{xiao2001prolate} and can be found in a slightly different form as
Theorem 4.1 in \cite{xiao2001prolate}.

\begin{theorem}
The $j$th smallest eigenvalue of $\mtx{A}_{\rm even}$ is $\chi_{n}$, where
$n=2j$.
Moreover, the corresponding eigenvector $\beta = \{ \beta_k \}_{k=0}^\infty$ satisfies
\begin{equation}
\alpha_{k} = \begin{cases} \beta_{k/2} \sqrt{k + \frac{1}{2}} & k \text{ is even} \\ 0 & k \text{
is odd} \end{cases}, \quad k = 0,1,\hdots
\label{eq:alphaeven}
\end{equation}
where $\{\alpha_k\}_{k=0}^\infty$ is given in \eqref{eq:legeexp}.
Similarly, the $j$th smallest eigenvalue of $\mtx{A}_{\rm odd}$ is
$\chi_{n}$, where $n=2j+1$, and 
the corresponding eigenvector $\beta = \{ \beta_k \}_{k=0}^\infty$ satisfies
\begin{equation}
\alpha_{k} = \begin{cases} 0 & k \text{ is even} \\ \beta_{(k-1)/2}\sqrt{k + \frac{1}{2}} & k
\text{ is odd}
\end{cases}, \quad k = 0,1,\hdots.
\label{eq:alphaodd}
\end{equation}
\label{thm:xiao}
\end{theorem}

\begin{observe}
In practice, the matrices $A_{\rm even}$ and $A_{\rm odd}$ must of course be
truncated before being diagonalized. 
We find that to compute $\psi_n$ for a given integer $n \geq 0$ to full
quadruple precision accuracy, it suffices to take the $M \times M$ leading
principal submatrix of the relevant matrix, where $M = \lceil 1.1c + n +
1000 \rceil$. 
\label{obs:truncate}
\end{observe}

\subsection{Computing eigenvalues of the truncated Fourier transform}

The following theorem provides a way of computing an eigenvalue of the
truncated Fourier transform, if the Legendre series expansion of the
respective PSWF is available.
It can be found as Theorem 3.27 in \cite{osipov2013numerical}. 

\begin{theorem}
Suppose that $\psi_n$ is a PSWF with Legendre expansion given by
\eqref{eq:legeexp}.
Then, 
\begin{equation}
\la_n = \frac{2 \alpha_0}{\psi_n(0)}
\label{eq:eig1}
\end{equation}
if $n$ is even, and
\begin{equation}
\la_n = \frac{2 c i \alpha_1}{3 \psi_n'(0)}
\label{eq:eig2}
\end{equation}
if $n$ is odd.
\label{thm:eig}
\end{theorem}

\begin{observe}
Theorem \ref{thm:eig} shows that computing $\lambda_n$ for large $n$
requires computing small entries of the eigenvector of a tridiagonal matrix.
At first glance, it seems that $\lambda_n$ can only be computed using
\eqref{eq:eig1} and \eqref{eq:eig2} when $|\lambda_n|$ is above the working
precision.
However, it is shown in \cite{osipov2017evaluation} that for a certain class
of symmetric tridiagonal matrices, which includes the matrices $A_{\rm
odd}$ and $A_{\rm even}$, the entries of the eigenvectors can be computed to
high relative accuracy using, for example, Rayleigh quotient iteration. 
This requires modifying the convergence criterion to track the convergence
of the small entries of the eigenvector, in addition to the eigenvalue.
This observation is crucial to Algorithm \ref{alg:eig} below.
\label{obs:osipov}
\end{observe}

\subsection{Computing eigenvalues of the differential operator}
\label{s:chi}
For any integer $n \geq 0$, $\chi_n$ can be computed using Theorem
\ref{thm:xiao}. 
We let $M = \lceil 1.1c + n + 1000 \rceil$ and compute the matrix $T \in
\R^{M \times M}$, consisting of the $M \times M$ leading principal submatrix
of $A_{\rm even}$, if $n$ is even, or $A_{\rm odd}$, if $n$ is odd.

We then perform Sturm bisection as described below.
First, we let $m = \lfloor n/2 \rfloor$, $a = 0$, and $b=(1+2n)c$.
After this, we compute $\rho(b)$ where $\rho$ is defined in Theorem
\ref{thm:sturm}, using Lemma \ref{lemma:charpol}.  
If $\rho(b) \leq m-1$, we set $b$ to twice its value, and again check if
$\rho(b) \leq m-1$.
This is repeated until $\rho(b) \geq m$.

We then set $d=(a+b)/2$ and evaluate $\rho(d)$.
If $\rho(d) \leq m-1$, we set $a$ to $d$, and if $\rho(d) \geq m$, we set
$b$ to $d$.
This procedure is continued until $\rho(a)=m-1$ and $\rho(b)=m$, and $b-a$
is below the desired accuracy.
At this point, we can take $\chi_n = (b+a)/2$.

\subsection{Root-finding for special functions}
\label{s:rootfinding}

We use a scheme introduced in \cite{glaser2007fast} due to its linear
complexity with respect to the number of roots, favorable constants, and
ability to attain high accuracy.
When an expansion of the form \eqref{eq:legeexp} is available, the scheme
determines all roots, one at a time, starting with an initial root. 

\noindent \textbf{Determining the first root}

\noindent
If $n$ is odd, then 0 is a root of $\psi_n$. 
We evaluate $\psi_n'(0)$ using the expansion \eqref{eq:diffexp}, and then
move onto determining subsequent roots.

If $n$ is even, then $\psi_n'(0) = 0$.
Thus, a crude approximation to the smallest positive root of $\psi_n$ can
be found by solving the ODE \eqref{eq:prufode2} subject to the initial
condition
\begin{equation}
x(0) = 0
\end{equation}
from 0 to $-\pi/2$.
In our implementation, we use a second-order Runge--Kutta scheme;
any other low-order, single-step scheme could be used. 
Let $\wt{x}^{(0)}$ be the approximation resulting from evaluating the
solution at $-\pi/2$. 
We can now refine it using the standard Newton iteration
\begin{equation}
 \wt{x}^{(k)} = \wt{x}^{(k-1)} -
 \frac{\psi_n(\wt{x}^{(k-1)})}{\psi_n'(\wt{x}^{(k-1)})}, \quad k =
 1,2,\hdots.
 \label{eq:newton}
\end{equation}
Naturally, this requires evaluating $\psi_n$ and $\psi_n'$ in the vicinity
of $\wt{x}^{(0)}$.
To do this, we first evaluate $\psi_n(0)$ using \eqref{eq:legeexp}, and then
determine higher-order derivatives using the recurrence \eqref{eq:derivs}.
After this, we use the Taylor series \eqref{eq:gen_taylor1} and
\eqref{eq:gen_taylor2} with $\psi=\psi_n$ and $x_* = 0$ to evaluate $\psi_n$
and $\psi_n'$ in \eqref{eq:newton}.
It was determined in \cite{glaser2007fast} that 30 and 60 terms in the
Taylor series suffice to get full double and quadruple precision,
respectively. 
After the Newton iteration has converged, the value of $\psi_n'$ at the
newly obtained root is determined using its Taylor series.

\noindent \textbf{Determining subsequent roots}

\noindent When a root $x_m$ is known, as well as the value of
$\psi_n'(x_m)$, the next largest root, $x_{m+1}$, is determined using a
two-step procedure.
In the first step, \eqref{eq:prufode2} is solved subject to the initial condition  
\begin{equation}
x \left( \frac{\pi}{2} \right) = x_m
\end{equation}
from $\pi/2$ to $-\pi/2$. 
Again, we use a second-order Runge--Kutta method, and $\wt{x}^{(0)}$, an
approximation to the true root, is obtained by evaluating the solution at
$-\pi/2$.
In the second step, higher-order derivatives at $x_m$ are obtained using
\eqref{eq:derivs}.
Then, the Newton iteration \eqref{eq:newton} is used to refine the
approximation to the root, where $\psi_n$ and $\psi_n'$ are evaluated using
\eqref{eq:gen_taylor1} and \eqref{eq:gen_taylor2} with $\psi=\psi_n$ and
$x_* = x_m$.
After $x_{m+1}$ is obtained, $\psi_n'(x_{m+1})$ is evaluated using
the Taylor series.

We repeat this procedure to determine all positive roots and then use
symmetry to determine negative roots.

\subsection{Quadratures for bandlimited functions}
\label{s:weights}

It turns out that the Euclidean algorithm for polynomial division has an
analog for bandlimited functions.
This is summarized by the following result, which is Theorem 6.2 in
\cite{xiao2001prolate}. 

\begin{theorem}
Suppose that $\sa, \phi : [-1,1] \to \C$ are a pair of $C^2$-functions such
that $\phi(1) \neq 0$ and $\phi(1) \neq 0$, and the functions $f$ and $p$
are defined by the formulas
\begin{align}
f(x) &= \int_{-1}^1 \sa(t) e^{2icxt}\,dt \\
p(x) &= \int_{-1}^1 \phi(t) e^{icxt}\,dt.
\end{align}
Then there exist two $C^1$-functions $\eta, \xi: [-1,1] \to \C$ such that
\begin{equation}
f(x) = p(x) q(x) + r(x), \quad x \in \R,
\end{equation}
with the functions $q,r :[0,1] \to \R$ defined by the formulas
\begin{align}
q(x) &= \int_{-1}^1 \eta(t) e^{icxt}\,dt \\
r(x) &= \int_{-1}^1 \xi(t) e^{icxt}\,dt.
\end{align}
\label{thm:euclid}
\end{theorem}
It follows from this result that for sufficiently large $n$, if the roots of
$\psi_n$ are chosen as quadrature nodes and the weights are chosen such that
the rule integrates $\{\psi_k\}_{k=0}^{n-1}$ exactly, the rule will
accurately integrate functions with bandlimit up to $2c$.
This is made rigorous by the following theorem which can be found in a
slightly different form Theorem 6.3 in \cite{xiao2001prolate}. 

\begin{theorem}
Suppose $n > 0$ is an integer, and 
let $\{x_j\}_{j=1}^n$ denote the roots of $\psi_n$. 
Suppose that the quadrature rule $\{(w_j,x_j)\}_{j=1}^n$ integrates exactly
the first $n$ PSWFs of bandlimit $c$, $\{\psi_k\}_{k=0}^{n-1}$. 
Then, we have the following bound for every function $f:[-1,1] \to \C$ that
satisfies the conditions of Theorem \ref{thm:euclid}: 
\begin{equation}
\left| \sum_{j=1}^n w_j f(x_j) - \int_{-1}^1 f(x)\,dx \right| \leq 
|\lambda_n| \| \eta \|_{2} + \|\xi\|_2 \sum_{k=n}^\infty |\lambda_k| \|
\psi_k\|_{\infty}^2 \left(2 + \sum_{j=1}^n |w_j| \right),
\end{equation}
where the functions $\eta, \xi: [-1,1] \to \C$ are defined in Theorem
\ref{thm:euclid} with $\phi=\psi_n$ and $\| \psi_k \|_\infty$ denotes the
$L^\infty([-1,1])$-norm of $\psi_k$. 
\label{thm:quad}
\end{theorem}

\subsection{Quadrature weights}
\label{s:rule}

Based on Theorem \ref{thm:quad}, we choose the roots $\{x_j\}_{j=1}^n$ of
$\psi_n$ to be the quadrature nodes.
We will order the roots so that $\{x_j\}_{j=1}^n$ is increasing. 
In this and the subsequent section, we describe a scheme 
for efficiently computing the corresponding quadrature weights
$\{w_j\}_{j=1}^n$ such that the rule $\{(w_j,x_j)\}_{j=1}^n$ integrates
$\{\psi_k\}_{k=0}^{n-1}$ accurately.
For this, we follow the approach in Section 9.4 of
\cite{osipov2013numerical}, which chooses the weights based on a
certain partial fractions approximation of the PSWFs.

We first define the functions 
\begin{equation} 
\varphi_j(x) = \frac{\psi_n(x)}{\psi_n'(x_j)(x-x_j)}, \quad j =
1,2,\hdots,n,
\label{eq:phi} 
\end{equation}
and then choose weights  
\begin{equation}
w_j = \int_{-1}^1 \varphi_j(x)\,dx, \quad j = 1,2,\hdots,n.
\label{eq:weights}
\end{equation}
We will use
\begin{equation}
\{(w_j, x_j)\}_{j=1}^n
\label{eq:rule}
\end{equation}
as a quadrature rule.

It turns out that for sufficiently large $n$, the rule
\eqref{eq:rule} will also integrate $\{\psi_k\}_{k=0}^{n-1}$ to high
accuracy.
This is summarized in the following theorem, which is found in a slightly
differnt form in Theorem 9.13 in \cite{osipov2013numerical}. 
\begin{theorem}
Let $c > 60$ and $\eps \in (0,1)$ be fixed, and let
$n > 0$ and $0 \leq k < n$  be integers, such that
\begin{equation}
n > \frac{2c}{\pi} + \left( 10 + \frac{3}{2} \log c + \frac{1}{2} \log
\frac{1}{\eps} \right) \log \left(\frac{c}{2}\right).
\end{equation}
Then
\begin{equation}
\left| \int_{-1}^1 \psi_k(x)\,dx - \sum_{j=1}^n w_j \psi_k(x_j) \right| <
\eps.
\end{equation}
\label{thm:weights}
\end{theorem}

\begin{remark}
We note that the choice of weights in \eqref{eq:weights} has a clear parallel
to Gauss--Legendre quadrature.  
If $\{x_j\}_{j=1}^n$ are the roots of $P_n$, the corresponding weights in
the Gauss--Legendre rule are chosen as the integrals of the Lagrange basis.
To be precise,
\begin{equation}
w_j = \int_{-1}^1 \ell_j(x)\,dx, \quad j = 1,2,\hdots,n,
\end{equation}
where
\begin{equation}
\ell_j(x) = \frac{P_n(x)}{P_n'(x_j)(x-x_j)}, \quad j = 1,2,\hdots,n.
\label{eq:lagrange}
\end{equation}
We note the similarities between \eqref{eq:lagrange} and \eqref{eq:phi}.
\end{remark}

\subsection{Computing the quadrature weights} 

It follows from Theorem \ref{thm:weights} that for sufficiently large $n$,
the rule \eqref{eq:rule} integrates $\{\psi_k\}_{k=0}^{n-1}$
virtually exactly. 
Let $\psi_n$ have a Legendre expansion given by \eqref{eq:legeexp}.
Using the identity \eqref{eq:legeid}, it follows that for $j =
1,2,\hdots,n$,
\begin{equation}
w_j = \int_{-1}^1 \varphi_j(x)\,dx = \frac{1}{\psi_n'(x_j)} \int_{-1}^1 \frac{\psi_j(x)}{x-x_j}
\,dx = -\frac{2}{\psi_n'(x_j)}\sum_{k=0}^\infty \alpha_k Q_k(x_j) =
-\frac{2 \Psi_n(x_j)}{\psi_n'(x_j)},
\end{equation}
where
\begin{equation}
\Psi_n(x) = \sum_{k=0}^\infty \alpha_k Q_k(x),
\label{eq:Psi}
\end{equation}
where $Q_k$ denotes the Legendre function of the second kind of index $k$
(see Section \ref{s:legendre}).
Since $\{\psi_n'(x_j)\}_{j=1}^n$ is determined during the root-finding,
computing the weights $\{w_j\}_{j=1}^n$ is essentially equivalent to
evaluating the series \eqref{eq:Psi} at
$x_1,x_2,\hdots,x_n$. 
An efficient mechanism for this is given in the following theorem, which 
can be found in a slightly different form as Theorem 9.15 in
\cite{osipov2013numerical}. 

\begin{theorem}
Let $\Psi_n$ be defined as in \eqref{eq:Psi}. Then $\Psi_n$ satisfies the ODE
\begin{equation}
(1-x^2) \Psi_n''(x) - 2x \Psi_n'(x) + (\chi_n - c^2 x^2) \Psi_n(x) = -c^2
\alpha_0 x - \frac{c^2 \alpha_1}{3}, \quad x \in (-1,1).
\label{eq:Psi_ode}
\end{equation}
\label{thm:Psi}
\end{theorem}
\begin{observe}
Theorem \ref{thm:Psi} offers a straightforward way to evaluate $\Psi$. 
We first note that repeated differentiation of \eqref{eq:Psi_ode} yields the
equations 
\begin{equation}
(1-x^2) \Psi_n'''(x) -4x \Psi_n''(x) + (\chi_n - c^2 x^2 - 2) \Psi_n'(x) -
2c^2 x \Psi_n(x) = -c^2 \alpha_0^{(n)}, \quad x \in (-1,1),
\label{eq:Psi_ode2}
\end{equation}
and for any integer $k \geq 2$
\begin{align}
(1-x^2) \Psi_n^{(k+2)}(x) &= 2x (1 + k) \Psi_j^{(k+1)}(x) + (k^2 + k + c^2
x^2 - \chi_n) \Psi_n^{(k)}(x) \nonumber \\
&+ 2c^2 k x \Psi_n^{(k-1)}(x) + c^2 k (k - 1) \Psi_n^{(k-2)}(x), \quad x \in (-1,1).
\label{eq:Psi_derivs}
\end{align}

We then let $m=\lfloor n/2+1 \rfloor$, so that $x_m$ is the smallest
non-negative root of $\psi_n$,
and determine the value of $\Psi_n(x_m)$ and $\Psi_n'(x_m)$ using
\eqref{eq:Psi}.
We then determine the values $\{\Psi_n(x_j)\}_{j=m+1}^{n-4}$ by solving
the ODE \eqref{eq:Psi_ode} using the Taylor series method described in
Section \ref{s:taylor} and the recurrence relations \eqref{eq:Psi_ode2} and
\eqref{eq:Psi_derivs}.
Again, 30 and 60 terms in the Taylor series suffice to get full accuracy in
double and quadruple precision, respectively. 
The Taylor series method loses accuracy near $x=1$.
As such, we evaluate $\{\Psi_n(x_j)\}_{j=n-3}^{n}$ directly using
\eqref{eq:Psi}.
The values $\{\Psi_n(x_j)\}_{j=1}^m$ are evaluated using symmetry. 

\label{obs:Psi}
\end{observe}

\section{Algorithm}
\label{s:algorithm}

\subsection{Overview of the algorithm} 

Given a bandlimit $c > 0$ and a sufficiently large integer $n > 0$, our
goal is to compute the $n$-point quadrature rule \eqref{eq:rule}. 
As intermediate steps, $\chi_n$ in \eqref{eq:ode} and $\{\alpha_j\}_{j=0}^N$
in \eqref{eq:legeexp} are computed. 
Moreover, the corresponding eigenvalue $\lambda_n$ of the operator $F_c$ in
\eqref{eq:fourop} is computed to high accuracy, which is useful since it
appears in various error bounds (e.g., see Theorems 9.6 and 9.8 in
\cite{osipov2013numerical}).  


{\bf Obtaining an approximation to $\mathbf{\chi_n}$.} 
We first compute $\wt{\chi}_n$ such that
\begin{equation}
    | \wt{\chi}_n - \chi_n| < \max_{m \neq n} | \wt{\chi}_n - \chi_m|.
    \label{eq:chibound}
\end{equation}
We let $M = \lceil 1.1 c + n + 1000\rceil$ (see Observation
\ref{obs:truncate}) and compute the leading $M \times M$ leading principal
submatrix of $\mtx{A}_{\rm even}$ using \eqref{eq:aeven1} and
\eqref{eq:aeven2} if $n$ is even and $\mtx{A}_{\rm odd}$ using
\eqref{eq:aodd1} and \eqref{eq:aodd2} if $n$ is odd.
Since these matrices are tridiagonal, only $O(n)$ entries need to be
computed. 
We then perform Sturm bisection as described in Section \ref{s:chi}.
Terminating the bisection when \eqref{eq:chibound} and $|\wt{\chi}_n -
\chi_n| < 2^{-40}$ works well in practice. 
While this requires a total number of operations that scales proportionately
to $n \log n$, the associated constant is very small.

{\bf Computing $\mathbf{\{ \alpha_j \}_{j=0}^N}$.}
Next, $\wt{\chi}_n$ is used as the initial shift for Rayleigh quotient
iteration on the leading $(2n+4) \times (2n+4)$ leading principal submatrix
of $\mtx{A}_{\rm odd}$ if $n$ is odd and $\mtx{A}_{\rm even}$ if $n$ is
even. 
The standard theory for Rayleigh quotient iteration (cf.~\cite[Theorem
27.3]{trefethen1997numerical}) states that convergence to the corresponding
eigenvector $\beta = \{\beta_j \}_{j=0}^{2n+3}$ occurs in $O(n)$ operations.
It is critical that the Rayleigh quotient iteration should be terminated
based on both convergence to the eigenvalue and small entries of the
eigenvector (see Observation \ref{obs:osipov} and Algorithm \ref{alg:eig}).

The coordinates of $\beta$ provide the coefficients
$\{\alpha_j\}_{j=0}^{4n+7}$ in the Legendre expansion \eqref{eq:legeexp},
and the corresponding eigenvalue is $\chi_n$ (see Theorem \ref{thm:xiao}).
We choose $N$ to be the smallest integer such that $|\alpha_k| < \eps_{\rm mach}$ for $k
> N$ where $\eps_{\rm mach}$ is machine precision and only keep the coefficients
$\{\alpha_j\}_{j=0}^N$. 

{\bf Determining $\mathbf{\{x_j\}_{j=1}^n}$ and $\mathbf{\{\psi_n'(x_j)\}_{j=1}^n}$.}
Once $\{\alpha_j\}_{j=0}^N$ is determined, the roots $\{x_j\}_{j=1}^n$ of
$\psi_n$ and the derivative of $\psi_n$ at the roots
$\{\psi_n'(x_j)\}_{j=1}^n$ are obtained via the algorithm of Section
\ref{s:rootfinding}.
As explained in Section \ref{s:rootfinding}, this entire procedure requires
$O(n)$ operations.

{\bf Computing $\mathbf{\{w_j\}_{j=1}^n}$.}
Once $\{x_j\}_{j=1}^n$ and $\{\psi_n'(x_j)\}_{j=1}^n$ are available, the
quadrature weights are determined via the scheme in Observation
\ref{obs:Psi}.
Again, this requires $O(n)$ operations.

{\bf Computing $\mathbf{\la_n}$.} 
Given $\{\alpha_j \}_{j=0}^N$, the eigenvalue $\la_n$ is evaluated
in $O(n)$ operations using Theorem \ref{thm:eig} with $\psi_n(0)$ and
$\psi_n'(0)$ evaluated using \eqref{eq:legeexp} and
\eqref{eq:diffexp}, respectively.
Per Observation \ref{obs:osipov}, the result is computed to high relative
accuracy. 

\subsection{Detailed description of the algorithm}

Algorithm \ref{alg:bisection} below computes an approximation $\wt{\chi}_n$
to $\chi_n$ using Sturm bisection.
This is then used in Algorithm \ref{alg:rayleigh}, which computes the true
value of $\chi_n$ and $\{\alpha_j\}_{j=0}^N$. 
These quantities are the input of Algorithm \ref{alg:roots}, which computes
$\{ x_j \}_{j=1}^n$ and $\{ \psi_n'(x_j) \}_{j=1}^n$.
Algorithm \ref{alg:weights} computes $\{w_j\}_{j=1}^n$, and Algorithm
\ref{alg:eig} is used to compute $\lambda_n$ to high relative accuracy.

\begin{algorithm}[Compute $\wt{\chi}_n$]~

\noindent Input: $c$, $n$ \\
\noindent Output: $\wt{\chi}_n$

\begin{enumerate}
\item[Step 1.] Let $M = \lceil 1.1c + n + 1000 \rceil$ and compute the entries
on the diagonal and subdiagonal of the $M \times M$ leading principal
submatrix of the matrix $\mtx{T}$, 
where $\mtx{T} = \mtx{A}_{\rm even}$ using \eqref{eq:aeven1} and
\eqref{eq:aeven2} if $n$ is even, and $\mtx{T} = \mtx{A}_{\rm odd}$ using
\eqref{eq:aodd1} and \eqref{eq:aodd2} if $n$ is odd. 

\item[Step 2.] Determine an upper bound on $\chi_n$ using the function
$\rho$, as defined in
Theorem \ref{thm:sturm}:

\begin{enumerate}

\item Set $a := 0$ and $b := (1+2n)c$, and let $m = \lfloor n/2 \rfloor$.

\item Evaluate $\rho(b)$ using Lemma \ref{lemma:charpol}. 

\item 
If $\rho(b) \geq m$, continue to Step 3. If not, set $a := b$ and
$b := 2b$ and go back to (a). 

\end{enumerate}

\item[Step 3.] Perform bisection on $\rho$:

\begin{enumerate}

\item Set $d := (a+b)/2$ and evaluate $\rho(d)$ using Lemma
\ref{lemma:charpol}. 

\item If $\rho(d) \leq m-1$, set $a := d$. If $\rho(d) \geq m$, then set $b
:= d$.

\item 
If converged, move to Step 4.
Otherwise, go back to (b).

\end{enumerate}

\item[Step 4.] Set $\wt{\chi}_n = (a+b)/2$.
\end{enumerate}

\label{alg:bisection}
\end{algorithm}

\begin{algorithm}[Compute $\{\alpha_j\}_{j=0}^N$]~

\noindent Input: $c$, $n$, $\wt{\chi}_n$ \\
\noindent Output: $\chi_n$, $\{ \alpha_j \}_{j=0}^N$ 

\begin{enumerate}
\item[Step 1.] Let $M = 2n + 4$ and compute the entries on the diagonal and
subdiagonal of the $M \times M$ leading principal submatrix of the matrix
$\mtx{T}$, where $\mtx{T} = \mtx{A}_{\rm even}$ using
\eqref{eq:aeven1} and \eqref{eq:aeven2} if $n$ is even, and $\mtx{T} =
\mtx{A}_{\rm odd}$ using \eqref{eq:aodd1} and \eqref{eq:aodd2} if $n$ is
odd. 

\item[Step 2.] Perform Rayleigh quotient iteration on $\mtx{T}$ using
$\wt{\chi}_n$ as the initial approximation to the eigenvalue:
\begin{enumerate}
\item Set $\la := \wt{\chi}_n$ and form a random, unit vector $\beta \in
\R^M$.
\item Compute $x := (\mtx{T} - \la \mtx{I})^{-1} \beta$, overwrite
$\beta := x / \| x \|$, and set
$\wt{\la} := \beta^T \mtx{T} \beta$.
\item If converged (see Observation \ref{obs:osipov}), set $\chi_n =
\wt{\lambda}$ and move to Step 3.  
If not, set $\la := \wt{\la}$, and go back to (b). 
\end{enumerate}

\item[Step 3.] Compute $\{\alpha_j\}_{j=0}^{2M-1}$ using \eqref{eq:alphaeven}
and \eqref{eq:alphaodd}, if $n$ is even and odd, respectively. 
Find the smallest $N$ such that $|\alpha_k| < \eps_{\rm mach}$ for all
$k > N$, where $\eps_{\rm mach}$ is machine precision, and
return $\chi_n$ and $\{\alpha_j\}_{j=0}^N$.
\end{enumerate}

\label{alg:rayleigh}
\end{algorithm}

\begin{algorithm}[Determine $\{x_j\}_{j=1}^n$ and
$\{\psi_n'(x_j)\}_{j=1}^n$]~

\noindent Input: $c$, $n$, $\{ \alpha_j \}_{j=0}^N$, $\chi_n$ \\
\noindent Output: $\{x_j\}_{j=1}^n$, $\{\psi_n'(x_j)\}_{j=1}^n$

\begin{enumerate}
\item[Step 1.] Set $m = \lfloor \frac{n}{2} \rfloor + 1$. 
If $n$ is odd, evaluate $\psi_n'(0)$ using \eqref{eq:diffexp}, and then set
$roots(m) = 0$ and $ders(m) = \psi_n'(0)$, and continue to Step 2. 
If $n$ is even, conduct the following procedure: \begin{enumerate}

\item Evaluate $\psi_n(0)$ using \eqref{eq:legeexp}.
Compute $\psi^{(k)}_n(0)$ for $k = 2,3,\hdots,K$ where $K = 30$ in double
precision and $K = 60$ in quadruple precision, using the recurrence
relation \eqref{eq:derivs} with the computed value of $\psi_n(0)$ and
$\psi_n'(0) = 0$.

\item Perform a crude solve on \eqref{eq:prufode2} subject to the initial
condition $x(0) = 0$ from $0$ to $-\pi/2$ using a second-order Runge--Kutta
scheme, and let $\wt{x}^{(0)}$ be the computed value of the solution at $\pi/2$. 

\item Refine $\wt{x}^{(0)}$ using the Newton iteration \eqref{eq:newton}.  At
each iteration, $\psi_n$ and $\psi_n'$ are evaluated using Taylor series
with derivatives computed in (a). 
Set the obtained root to $roots(m)$.

\item Set $ders(m)$ to the value of $\psi_n'(roots(m))$, evaluated using the
Taylor series.

\end{enumerate}

\item[Step 2.] The remaining non-negative roots can be found through the
following procedure:
\begin{enumerate}
\item Set $j := m$.
\item Evaluate $\psi^{(k)}_n(roots(j))$ for $k = 2,3,\hdots,K$, where $K=30$
in double precision and $K=60$ in quadruple precision, using the recurrence
relation \eqref{eq:derivs} and the facts that $\psi_n(roots(j)) = 0$ and
$\psi_n'(roots(j)) = ders(j)$.

\item
Perform a crude solve on \eqref{eq:prufode2} subject to the initial
condition $x(\pi/2) = roots(j)$ from $\pi/2$ to $-\pi/2$ using a
second-order Runge--Kutta scheme, and let $\wt{x}^{(0)}$ be the value of the
resulting solution at $-\pi/2$.

\item Refine $\wt{x}^{(0)}$ using the Newton iteration \eqref{eq:newton}, where
$\psi_n$ and $\psi_n'$ are evaluated by their Taylor series expansions with
derivatives computed in (b). 
Set the resulting root to $roots(j+1)$.

\item Determine $ders(j+1)$ by evaluating $\psi_n'(roots(j+1))$ using the
Taylor series.

\item If $j = n - 1$, continue to Step 3. Otherwise, set $j := j + 1$ and
return back to (b). 

\end{enumerate}

\item[Step 3.] For $j=1,2,\hdots,m-1$, set $roots(j) = -roots(n-j+1)$ and $ders(j) = (-1)^{n+1}
ders(n-j+1)$.

\item[Step 4.] Set $x_j = roots(j)$ and note that $\psi'_n(x_j) = ders(j)$
for $j = 1,2,\hdots,n$.

\end{enumerate}

\label{alg:roots}
\end{algorithm}

\begin{algorithm}[Compute $\{w_j\}_{j=1}^n$]~

\noindent Input: $c$, $n$, $\{\alpha_j\}_{j=0}^N$, $\chi_n$, $\{x_j\}_{j=1}^n$,
$\{\psi_n'(x_j)\}_{j=1}^n$ \\
\noindent Output: $\{w_j\}_{j=1}^n$

\begin{enumerate}
\item[Step 1.] Set $m = \lfloor \frac{n}{2} \rfloor + 1$. Evaluate
$\Psi_n(x_m)$ and $\Psi_n'(x_m)$ using \eqref{eq:Psi}, and set these to
$vals(m)$, respectively.

\item[Step 2.] Next, determine the values of $\Psi_n(x_j)$ and $\Psi_n'(x_j)$
for $j = m+1,m+2,\hdots,n-4$ using the following procedure:

\begin{enumerate}
\item Set $j := m$.
\item Evaluate $\Psi_n^{(k)}(x_j)$ for $k = 3,4,\hdots,K$, where $K = 30$ in
double precision and $K = 60$ in quadruple precision, using
\eqref{eq:Psi_ode}, \eqref{eq:Psi_ode2}, \eqref{eq:Psi_derivs}, and
 the initial conditions $\Psi_n(x_j) = vals(j)$ and $\Psi_n'(x_j) =
 ders(j)$.
\item Use the derivatives computed in (b) to compute Taylor series for
$\Psi_n(x_{j+1})$ and $\Psi_n'(x_{j+1})$ and set these to $vals(j+1)$ and
$ders(j+1)$, respectively. 
\item If $j = n-5$, continue to Step 3. If not, set $j := j + 1$ and go back to
(b). 
\end{enumerate}
\item For $j = n-3,n-2,n-1,n$, obtain $vals(j)$ by directly evaluating
$\Psi_n(x_j)$ using the series \eqref{eq:Psi}.

\item Set $vals(j) = (-1)^{n+1} vals(n-j+1)$ for $j = 1,2,\hdots,m-1$.

\item Finally, set $w_j = -2 vals(j)/\psi_n'(x_j)$ for $j = 1,2,\hdots,n$.

\end{enumerate}

\label{alg:weights}
\end{algorithm}

\begin{algorithm}[Compute $\la_n$]~

\noindent Input: $c$, $n$, $\{\alpha_j\}_{j=0}^N$\\
\noindent Output: $\la_n$

\begin{enumerate}
\item[Step 1.] If $n$ is even, evaluate $\psi_n(0)$ using \eqref{eq:legeexp}, and let $\la_n = 2\alpha_0 /
\psi_n(0)$. If $n$ is odd, evaluate $\psi_n'(0)$ using 
\eqref{eq:diffexp} and let $\la_n = 2ci \alpha_1 / (3 \psi_n'(0))$.
\end{enumerate}

\label{alg:eig}
\end{algorithm}

\section{Numerical experiments}

The algorithm described in Section \ref{s:algorithm} for computing the
quadrature rule \eqref{eq:rule} was implemented in Fortran with the version 
11.2.0 GNU Fortran compiler.  
The following experiments were run on a 2.1 GHz laptop with 16 GB of RAM.

\label{s:numerics}

\subsection{Experiment 1}
\label{s:exp1}
In this experiment, we measure the accuracy of and time to compute the
quadrature rule \eqref{eq:rule} as functions of the bandlimit $c$ and the
requested accuracy $\eps$. 
For any $c > 0$, we define the points 
\begin{equation}
\omega^c_k = \frac{2kc}{100}, \quad  k = 1,2,\hdots,100.
\label{eq:omegadef}
\end{equation}
For any integer $n \geq 0$, we define the error 
\begin{equation}
E_n^c = \max_{k=1,2,\hdots,100} \left| \int_{-1}^1 \cos(\omega_k^c x)\,dx -
\sum_{k=1}^n w_j^{(n)} \cos(\omega_k^c x_j^{(n)}) \right|,
\label{eq:err1}
\end{equation}
where $\{(w_j^{(n)},x_j^{(n)})\}_{j=1}^n$ is the $n$-point rule
\eqref{eq:rule}. 
We also define $n(\eps)$ to be smallest integer $n$ such that 
\begin{equation}
|\lambda_n| < \eps,
\label{eq:neps}
\end{equation}
for any $\eps \in (0,1)$.

For fixed values of $c$ and $\varepsilon$, we compute the
$n(\varepsilon)$-point quadrature rule as described in Section
\ref{s:algorithm} and then check $E_{n(\varepsilon)}^c$.
The results are shown in Table \ref{tab:exp1}, which has the following the structure. 
The first and second columns show the values of $c$ and $\varepsilon$, respectively.
The third column shows $n(\eps)$ as defined by \eqref{eq:neps}.
The fourth column shows $|\lambda_{n(\eps)}|$, computed in double precision.
The fifth and sixth columns show the values of $E_{n(\varepsilon)}^c$ when
the quadrature rule is computed in double and quadruple precision,
respectively. 
The seventh and eighth columns show the total CPU time taken in seconds
$t_{\rm total}$ to compute the quadrature rule in double and quadruple
precision, respectively.

The results show that for all choices of parameters such that $\eps > c
\eps_{\rm mach}$, where $\eps_{\rm mach}$ is machine precision,
$E_{n(\eps)}^c < \eps$. 
When $\eps < c\eps_{\rm mach}$, $E_{n(\eps)}^c$ is limited to approximately $c
\eps_{\rm mach}$. 
This is expected since the evaluation of the integrand in \eqref{eq:err1}
has condition number approximately $c$.
Thus, the quadrature rule is capable of integrating functions down to
machine precision.
Moreover, $|\lambda_n|$ provides a good upper bound for the error when the
$n$-point rule is used (see page 344 in \cite{osipov2013numerical} for
further discussion).  

\begin{table}[t]
\centering
\begin{tabular}{c|c|c|c|c|c|c|c}
$c$ & $\varepsilon$ & $n(\varepsilon)$ & $|\lambda_{n(\varepsilon)}|$ & $E_{n(\varepsilon)}^c$ (double) & $E_{n(\varepsilon)}^c$ (quad) & $t_{\rm total}$ (double) & $t_{\rm total}$ (quad) \\
\hline
$10^2$ & $10^{-10}$ & 86 & 0.59988E-10 & 0.49E-12 & 0.49E-12 & 0.554E-03 & 0.265E-01\\
$10^2$ & $10^{-25}$ & 112 & 0.33640E-25 & 0.28E-14 & 0.56E-28 & 0.653E-03 & 0.290E-01 \\
$10^2$ & $10^{-50}$ & 147 & 0.44641E-50 & 0.14E-13 & 0.66E-32 & 0.739E-03 & 0.332E-01
\\
\hline
$10^3$ & $10^{-10}$ & 667 & 0.95582E-10 & 0.27E-11 & 0.27E-11 & 0.238E-02 & 0.102E+00\\
$10^3$ & $10^{-25}$ & 708 & 0.97844E-25 & 0.24E-13 & 0.32E-28 & 0.239E-02 & 0.109E+00\\
$10^3$ & $10^{-50}$ & 768 & 0.39772E-50 & 0.12E-13 & 0.81E-32 & 0.261E-02 & 0.117E+00
\\
\hline 
$10^4$ & $10^{-10}$ & 6405 & 0.57608E-10 & 0.35E-12 & 0.18E-12 & 0.203E-01 & 0.854E+00\\
$10^4$ & $10^{-25}$ & 6462 & 0.63792E-25 & 0.42E-12 & 0.68E-29 & 0.204E-01 & 0.859E+00\\
$10^4$ & $10^{-50}$ & 6548 & 0.51349E-50 & 0.15E-12 & 0.16E-30 & 0.200E-01 & 0.877E+00\\
\hline
$10^5$ & $10^{-10}$ & 63707 & 0.71063E-10 & 0.83E-11 & 0.33E-13 & 0.188E+00 & 0.884E+01 \\
$10^5$ & $10^{-25}$ & 63780 & 0.92981E-25 & 0.11E-10 & 0.63E-29 & 0.191E+00 & 0.887E+01 \\
$10^5$ & $10^{-50}$ & 63893 & 0.80840E-50 & 0.44E-11 & 0.21E-29 & 0.210E+00 & 0.933E+01\\
\hline
$10^6$ & $10^{-10}$ & 636670 & 0.79326E-10 & 0.19E-08 & 0.15E-12 & 0.184E+01 & 0.855E+02\\
$10^6$ & $10^{-25}$ & 636760 & 0.77413E-25 & 0.43E-09 & 0.83E-29 & 0.184E+01 & 0.855E+02\\
$10^6$ & $10^{-50}$ & 636900 & 0.69235E-50 & 0.29E-10 & 0.14E-27 & 0.187E+01 & 0.858E+02\\
\hline
$10^7$ & $10^{-10}$ & 6366252 & 0.87469E-10 & 0.42E-08 & 0.16E-11 & 0.205E+02 & 0.939E+03 \\
$10^7$ & $10^{-25}$ & 6366358 & 0.97995E-25 & 0.20E-08 & 0.11E-25 & 0.225E+02 & 0.100E+04\\
$10^7$ & $10^{-50}$ & 6366525 & 0.91559E-50 & 0.83E-10 & 0.61E-27 & 0.223E+02 & 0.945E+03
\end{tabular}
\caption{Experiment 1}
\label{tab:exp1}
\end{table}

\subsection{Experiment 2}
\label{s:exp2}

\begin{figure}[t]
\begin{tabular}{cc}
\includegraphics{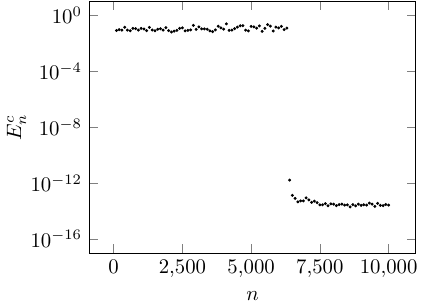} &
\includegraphics{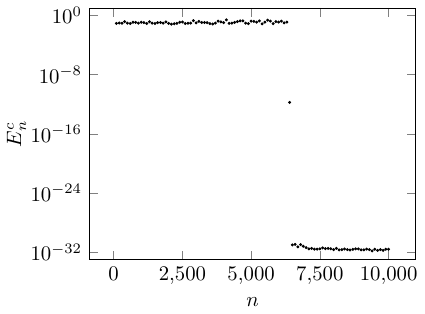}
\end{tabular}
\caption{$E^c_n$ in Experiment 2 when computations are done in double (left)
and quadruple (right) precision}
\label{fig:exp2_1}
\end{figure}

We now fix $c=10^4$ and examine how $E_n^c$ as defined in \eqref{eq:err1}
and $|\lambda_n|$ change as functions of $n$.
Figure \ref{fig:exp2_1} plots the values of $E_n^c$ when the quadrature rule
\eqref{eq:rule} is computed in double (left) and quadruple (right)
precision.
As is expected from Theorem \ref{thm:weights}, $E_n^c$ starts to decay
rapidly when $n \approx 2c/\pi$ and remains at approximately $c \eps_{\rm
mach}$ from then on. 

\begin{figure}[t]
\begin{tabular}{cc}
\includegraphics{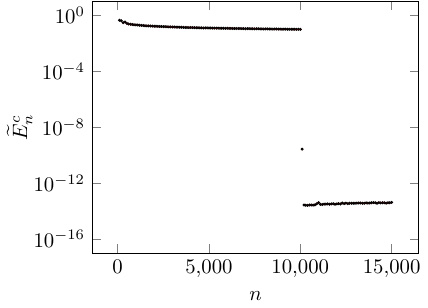} &
\includegraphics{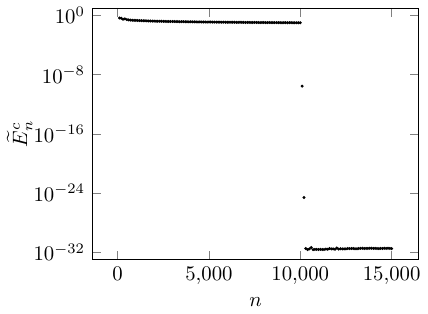}
\end{tabular}
\caption{$\widetilde{E}^c_n$ in Experiment 2 when computations are done in
double (left) and quadruple (right) precision}
\label{fig:exp2_3}
\end{figure}

For comparison, we define
\begin{equation}
\widetilde{E}_n^c = \max_{k=1,2,\hdots,100} \left| \int_{-1}^1 \cos(\omega_k^c x)\,dx -
\sum_{k=1}^n \widetilde{w}_j^{(n)} \cos(\omega_k^c \widetilde{x}_j^{(n)}) \right|,
\label{eq:err_gl}
\end{equation}
where $\{(\widetilde{w}_j^{(n)},\widetilde{x}^{(n)}_j)\}_{j=1}^n$ is the
$n$-point Gauss--Legendre quadrature rule and $\omega_k^c$ is defined in \eqref{eq:omegadef}.
Figure \ref{fig:exp2_3} plots the values of $\widetilde{E}_n^c$ when
$c=10^4$ in double (left) and quadruple (right) precision.
$\widetilde{E}_n^c$ starts to decay when $n \approx c$, which is a factor of
$\pi/2$ larger than than the point where $E_n^c$ starts to decay.

\begin{figure}[t]
\begin{tabular}{cc}
\includegraphics{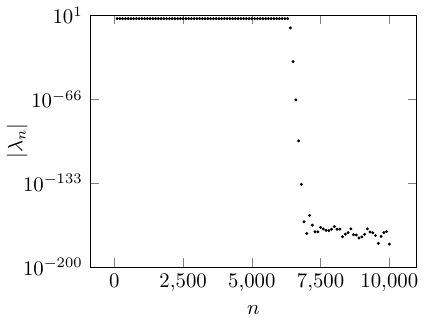} &
\includegraphics{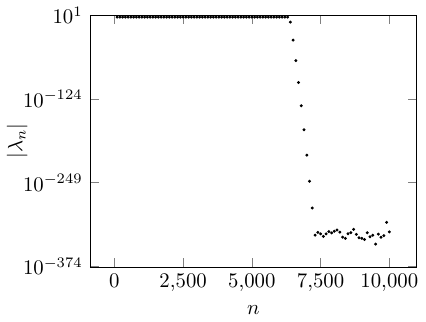}
\end{tabular}
\caption{$|\lambda_n|$ for Experiment 2 when computations are done in double
(left) and quadruple (right) precision}
\label{fig:exp2_2}
\end{figure}

Figure \ref{fig:exp2_2} plots the values of $|\lambda_n|$ for
varying values of $n$ in double (left) and quadruple (right) precision.
These values also start to decay superalgebraically when $n \approx 2c/
\pi$. 
As discussed in Observation \ref{obs:osipov}, we are able to compute
$|\lambda_n|$ well below machine precision, down to approximately the square
root of the point of underflow.

\subsection{Experiment 3}

\begin{figure}[t]
\begin{tabular}{cc}
\includegraphics{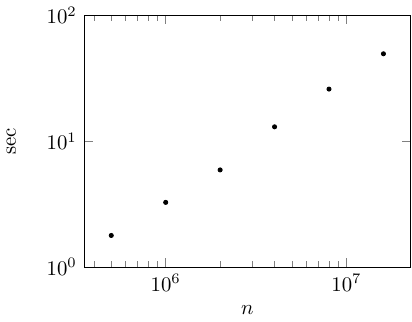} & 
\includegraphics{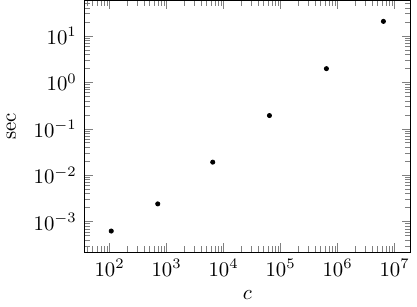}
\end{tabular}
\caption{$t_{\rm total}$ for computations done in double precision in
Experiments 3 (left) and 4 (right)}
\label{fig:exp4}
\end{figure}

In this experiment, we fix $c = 10^6$ and examine the CPU time required to
compute the $n$-point quadrature rule \eqref{eq:rule} as a function of $n$.
The results are displayed in Tables \ref{tab:exp3} and \ref{tab:exp3_ext}
when the computations are done in double and quadruple precision,
respectively.
Both tables have the following structure.
The first column contains the values of $n$. 
The second column shows $t_{\rm prol}$, which is the time required to compute
the Legendre expansion for $\psi_n$ using Algorithms \ref{alg:bisection} and
\ref{alg:rayleigh}.
The third column contains $t_{\rm roots}$, which is the time required to
determine the quadrature nodes using Algorithm \ref{alg:roots}.
The fourth column contains $t_{\rm weights}$, which is the time required to
compute the quadrature weights via Algorithm \ref{alg:weights}. 
Finally, the fifth column reports $t_{\rm total}$, which is the total time
taken to compute the quadrature rule. 
The left plot in Figure \ref{fig:exp4} shows the values of $t_{\rm total}$
plotted as a function of $n$ when the computation is done in double
precision.
We see that $t_{\rm total}$ scales approximately linearly with $n$, as
expected. 
The constants are significantly higher when the computation
is done in quadruple precision, which is expected since the machine these
experiments were run on does not support quadruple precision floating point
operations in hardware. 

\begin{table}[t]
\centering
\begin{tabular}{c|c|c|c|c}
$n$ & $t_{\rm prol}$ & $t_{\rm roots}$ & $t_{\rm weights}$ &
$t_{\rm total}$ \\
\hline
500,000 & 0.763E+00 & 0.834E+00 & 0.183E+00 & 0.178E+01 \\
1,000,000 & 0.124E+01 & 0.168E+01 & 0.337E+00 & 0.326E+01 \\ 
2,000,000 & 0.184E+01 & 0.338E+01 & 0.689E+00 & 0.591E+01 \\
4,000,000 & 0.404E+01 & 0.742E+01 & 0.154E+01 & 0.130E+02 \\
8,000,000 & 0.680E+01 & 0.160E+02 & 0.318E+01 & 0.259E+02 \\
16,000,000 & 0.136E+02 & 0.296E+02 & 0.628E+01 & 0.495E+02 
\end{tabular}
\caption{Experiment 3 (double precision)}
\label{tab:exp3}
\end{table}

\begin{table}[t]
\centering
\begin{tabular}{c|c|c|c|c}
$n$ & $t_{\rm prol}$ & $t_{\rm roots}$ & $t_{\rm weights}$ &
$t_{\rm total}$ \\
\hline
500,000 &  0.336E+02 & 0.559E+02 & 0.110E+02 & 0.101E+03 \\
1,000,000 &  0.453E+02 & 0.114E+03 & 0.230E+02 & 0.182E+03\\ 
2,000,000 &  0.545E+02 & 0.175E+03 & 0.338E+02 & 0.264E+03 \\
4,000,000 &  0.890E+02 & 0.402E+03 & 0.774E+02 & 0.569E+03 \\
8,000,000 &  0.191E+03 & 0.828E+03 & 0.155E+03 & 0.117E+04 \\
16,000,000 & 0.375E+03 & 0.165E+04 & 0.308E+03 & 0.234E+04 
\end{tabular}
\caption{Experiment 3 (quadruple precision)}
\label{tab:exp3_ext}
\end{table}

\subsection{Experiment 4}
We now fix the desired accuracy $\eps=e^{-50}$ and measure the CPU time
required to compute the quadrature rule \eqref{eq:rule} with $n(\eps)$
points (see \eqref{eq:neps}) as a function of the bandlimit $c$.
The results are displayed in Tables \ref{tab:exp4} and \ref{tab:exp4_ext}
when the computations are done in double and quadruple precision,
respectively.
Both tables have the following structure.
The first and second columns show the values of $c$ and $n(\eps)$,
respectively. 
The third column contains $t_{\rm prol}$, which is the time required to compute
the Legendre expansion for $\psi_n$ using Algorithms \ref{alg:bisection} and
\ref{alg:rayleigh}.
The fourth column contains $t_{\rm roots}$, which is the time required to
determine the quadrature nodes using Algorithm \ref{alg:roots}.
The fifth column contains $t_{\rm weights}$, which is the time required to
compute the quadrature weights via Algorithm \ref{alg:weights}. 
Finally, the sixth column reports $t_{\rm total}$, which is the total time
taken to compute the quadrature rule. 
The right plot in Figure \ref{fig:exp4} shows the values of $t_{\rm total}$
plotted as a function of $c$ when the computation is done in double
precision.
The data show that $t_{\rm total}$ asymptotically scales approximately
linearly with $c$. 

\begin{table}[t]
\centering
\begin{tabular}{c|c|c|c|c|c}
$c$ & $n(\eps)$ & $t_{\rm prol}$ & $t_{\rm roots}$ & $t_{\rm weights}$ &
$t_{\rm total}$ \\
\hline
$10^2$ & 107 & 0.359E-03 & 0.222E-03 & 0.440E-04 & 0.625E-03 \\
$10^3$ & 700 & 0.741E-03 & 0.138E-02 & 0.280E-03 & 0.240E-02 \\
$10^4$ & 6450 &  0.437E-02 & 0.121E-01 & 0.256E-02 & 0.191E-01 \\
$10^5$ & 63765 & 0.456E-01 & 0.115E+00 & 0.311E-01 & 0.192E+00 \\
$10^6$ & 636741 & 0.640E+00 & 0.110E+01 & 0.224E+00 & 0.197E+01 \\
$10^7$ & 6366336 & 0.481E+01 & 0.134E+02 & 0.227E+01 & 0.205E+02
\end{tabular}
\caption{Experiment 4 (double precision)}
\label{tab:exp4}
\end{table}

\begin{table}[t]
\centering
\begin{tabular}{c|c|c|c|c|c}
$c$ & $n(\eps)$ & $t_{\rm prol}$ & $t_{\rm roots}$ & $t_{\rm weights}$ &
$t_{\rm total}$ \\
\hline
$10^2$ & 107 & 0.165E-01 & 0.113E-01 & 0.244E-02 & 0.303E-01 \\
$10^3$ & 700 & 0.427E-01 & 0.653E-01 & 0.130E-01 & 0.121E+00 \\
$10^4$ & 6450 & 0.295E+00 & 0.632E+00 & 0.132E+00 & 0.106E+01  \\
$10^5$ & 63765 & 0.263E+01 & 0.901E+01 & 0.170E+01 & 0.133E+02 \\
$10^6$ & 636741 & 0.377E+02 & 0.710E+02 & 0.124E+02 & 0.121E+03 \\
$10^7$ & 6366336 & 0.371E+03 & 0.663E+03 & 0.129E+03 & 0.116E+04
\end{tabular}
\caption{Experiment 4 (quadruple precision)}
\label{tab:exp4_ext}
\end{table}

\section{Conclusions}
\label{s:conclusions}

In this paper, we describe a fast algorithm for computing quadrature
rules for bandlimited functions. 
For any fixed $c > 0$, the algorithm computes an $n$-point quadrature rule
in $O(n \log n)$ operations.
Moreover, the computed quadrature rules are capable of achieving machine
precision. 
The accuracy and speed of the algorithm are demonstrated in several
numerical experiments.
The approach of this paper can be generalized to design algorithms for
computing quadrature rules for a broad class of other special functions. 

\bibliographystyle{plain}
\bibliography{refs} 

\begin{thebibliography}{10}

\bibitem{abramowitz1988handbook}
Milton Abramowitz, Irene~A Stegun, and Robert~H Romer.
\newblock Handbook of mathematical functions with formulas, graphs, and
  mathematical tables, 1988.

\bibitem{boyd2004prolate}
John~P Boyd.
\newblock Prolate spheroidal wavefunctions as an alternative to {C}hebyshev and
  {L}egendre polynomials for spectral element and pseudospectral algorithms.
\newblock {\em Journal of Computational Physics}, 199(2):688--716, 2004.

\bibitem{glaser2007fast}
Andreas Glaser, Xiangtao Liu, and Vladimir Rokhlin.
\newblock A fast algorithm for the calculation of the roots of special
  functions.
\newblock {\em SIAM Journal on Scientific Computing}, 29(4):1420--1438, 2007.

\bibitem{landau1961prolate}
Henry~J Landau and Henry~O Pollak.
\newblock Prolate spheroidal wave functions, {F}ourier analysis and uncertainty
  {II}.
\newblock {\em Bell System Technical Journal}, 40(1):65--84, 1961.

\bibitem{landau1962prolate}
Henry~J Landau and Henry~O Pollak.
\newblock Prolate spheroidal wave functions, {F}ourier analysis and uncertainty
  {III}.
\newblock {\em Bell System Technical Journal}, 41(4):1295--1336, 1962.

\bibitem{osipov2017evaluation}
Andrei Osipov.
\newblock Evaluation of small elements of the eigenvectors of certain symmetric
  tridiagonal matrices with high relative accuracy.
\newblock {\em Applied and Computational Harmonic Analysis}, 43(2):173--211,
  2017.

\bibitem{osipov2013numerical}
Andrei Osipov, Vladimir Rokhlin, and Hong Xiao.
\newblock {\em Prolate Spheroidal Wave Functions of Order Zero}.
\newblock Springer, 2013.

\bibitem{shkolnisky2006approximation}
Yoel Shkolnisky, Mark Tygert, and Vladimir Rokhlin.
\newblock Approximation of bandlimited functions.
\newblock {\em Applied and Computational Harmonic Analysis}, 21(3):413--420,
  2006.

\bibitem{slepian1964prolate}
David Slepian.
\newblock Prolate spheroidal wave functions, {F}ourier analysis and uncertainty
  {IV}.
\newblock {\em Bell System Technical Journal}, 43(6):3009--3057, 1964.

\bibitem{slepian1978prolate}
David Slepian.
\newblock Prolate spheroidal wave functions, {F}ourier analysis, and
  uncertainty {V}.
\newblock {\em Bell System Technical Journal}, 57(5):1371--1430, 1978.

\bibitem{slepian1961prolate}
David Slepian and Henry~O Pollak.
\newblock Prolate spheroidal wave functions, {F}ourier analysis and uncertainty
  {I}.
\newblock {\em Bell System Technical Journal}, 40(1):43--63, 1961.

\bibitem{trefethen1997numerical}
Lloyd~N Trefethen and David Bau.
\newblock {\em Numerical linear algebra}, volume~50.
\newblock SIAM, 1997.

\bibitem{xiao2001prolate}
Hong Xiao, Vladimir Rokhlin, and Norman Yarvin.
\newblock Prolate spheroidal wavefunctions, quadrature and interpolation.
\newblock {\em Inverse problems}, 17(4):805, 2001.

\end{thebibliography}

\end{document}